\documentclass[preprint,authoryear,11pt]{elsarticle}
\usepackage{amsfonts,amssymb,amsmath,amsopn,amsthm}
\usepackage{graphics}
\usepackage{rotating}
\usepackage{bigints}
\usepackage{booktabs}
\usepackage{subfigure}
\usepackage{url}
\usepackage{lscape}
\usepackage{amssymb}
\usepackage{color}
\usepackage{textgreek}
\usepackage{enumerate}

\usepackage[colorlinks=true]{hyperref}
\usepackage{setspace}

\newcommand{\ignore}[1]{}
\newtheorem{asu}{{\sc Assumption}}

\newtheorem{theorem}{\sc Theorem}
\newtheorem{remark}{\sc Remark}
\newtheorem{proposition}{\sc Proposition}

\newdefinition{rmk}{Remark}

\usepackage[margin=1in]{geometry}
\newproof{pf}{Proof}
\newproof{pot}{Proof of Theorem \ref{thm2}}

\usepackage{amssymb}
\journal{Journal of Theoretical Probability}
\begin{document}
\begin{frontmatter}
\title{Large Deviation principles of Realized Laplace Transform of Volatility}
%\ignore{
\author[SDU]{Xinwei Feng}
\author[UM]{Lidan He\corref{cor1}}
\author[UM,ZUMRI]{Zhi Liu}
\address[SDU]{Zhongtai Securities Institute for Financial Studies, Shandong University, Jinan, Shandong Province, China}
\address[UM]{University of Macau, Macau SAR, China}
\address[ZUMRI]{Zhuhai-UM Science and Technology Research Institute, Zhuhai, China}
\cortext[cor1]{Corresponding author. Email: liuzhi@umac.mo.}
\begin{abstract}
Under scenario of high frequency data, consistent estimator of realized Laplace transform of volatility is proposed by \citet{TT2012a} and related central limit theorem has been well established. In this paper, we investigate the asymptotic tail behaviour of the empirical realized Laplace transform of volatility (ERLTV). We establish both large deviation principle and moderate deviation principle for the ERLTV. The good rate function for the large deviation principle is well defined in the whole real space, which indicates a limit for the normalized logarithmic tail probability of the ERLTV. Moreover, we also derive the function-level large and moderate deviation principles for ERLTV.
\end{abstract}
\end{frontmatter}
%\begin{spacing}{1.4}
\section{Introduction}
For a one-dimensional It$\hat{\text{o}}$ semi-martingale of the form
\begin{eqnarray}\label{model1}
    \mathrm{d}X_t=a_t \mathrm{d}t+\sigma_t \mathrm{d}W_t+\int_{\mathbb R}\delta(t-, x)\mu( \mathrm{d}t, \mathrm{d}x),\ t\geq 0,
\end{eqnarray}
on a probability space (the details will be specified later), one of quantities of interest for the process $\{X_t\}_{t\geq 0}$ is the quadratic variation of the continuous part, namely,  $\int_0^T\sigma_t^2\mathrm{d}t$ and $\{\sigma_t\}_{t\geq 0}$ itself is called volatility process. The quadratic variation which is called integrated volatility in the continuous-time modelling of financial markets, plays a key role in risk management, option pricing, and in high frequency trading. It has stimulated wide studies of statistical inference for the quadratic variation. The widespread use of continuous-time processes in macroeconomics and finance shows that time-varying volatility process is a salient empirical feature of many economic and financial time series. Because the underlying volatility process is latent and in most applications the available data are discrete, the statistical inference of volatility is complicated which can be uniquely determined by the discretized observations. Fortunately, the availability of high-frequency data makes it possible that using the  realized volatility to approximate the actual volatility. For details, see \citet{BNS2002}, \citet{ABM2005}, \citet{ABFO2010}, \citet{J2008}, \citet{MZ2009, MZ2016}, among others.

The quadratic variation characteristics the second order property of the process $X$, itself can be a stochastic process too. Recently, the property of quadratic variation attracts wide research interests. For instance, \citet{ABDE2001} investigated the distribution of volatility, \citet{V2015} studied the inference of the integrated volatility of volatility. In this paper, we are interested in the realized Laplace transform of volatility (RLTV) proposed in \citet{TT2012a}, where
$$RLTV_t=\int_0^t\exp(-u\sigma^2_s)\mathrm{d}s,~~~u\in \mathbb R^+.$$
RLTV contains more information than the integrated volatility (IV), because RLTV is a mapping from the real high frequency data to a random function (in terms of $u$) in $\mathbb R^+$ whereas the IV only maps the data to a random variable. Moreover, the further average of RLTV's over long time span preserves the complete information about the characteristics of volatility under some assumptions of stationarity, this provides a possibility to infer the distributional properties of the volatility. The statistical inference of RLTV is firstly considered by \citet{TTG2011} under a jump diffusion framework and the formal theory has been established by \citet{TT2012a}. The extension of RLTV to the framework of pure jump semi-martingales has been considered by \citet{TT2012b}, and the effect of microstructure noise has been studied by \citet{WLX2019} and \citet{WLX20192}. The bootstrap of the RLTV has been studies by \citet{HLV2020}.

Assuming that $\{X_t\}_{t\geq 0}$ is observed at some discrete points $\{t_i=i\Delta_n, i=0,1,\cdots\}$, \citet{TT2012a} proposed a consistent estimator (ERLTV) of RLTV. The related functional central limit theorem for the random function sequence has been established.

The central limit theorem, though is useful in approximating the finite sample distribution particularly around the center, it is less precise when estimating the tail probability. The purpose of this paper is to furnish some further estimations about the tail behaviour of sample distribution of the ERLTV, refining the functional central limit theorem obtained in \citet{TT2012a}. We establish the large and moderate deviation principles of ERLTV that provides a more accurate estimate of the tail probability of the ERLTV.

The large and moderate deviations problems arise in the theory of statistical inference, whose results provide the rates of convergence hence is useful in constructing asymptotic confidence intervals. See \citet{BM1993} for a comprehensive introduction of the application of large deviations in the statistical test theory. The large derivation and moderate deviation principles for the diffusion processes have been studied in existing literature. Limit to our knowledge, \citet{M2008} derived the large deviation for the threshold estimator for the constant volatility; \citet{DGW1999} obtained the large and moderate deviations for the quadratic variational processes of diffusions assuming the volatility process is time varying but non-random; \citet{KO2011} established the conditional large and moderate deviation principles of the realized volatility conditioning the path of the volatility; \citet{DS2014} and \citet{DGS2017} obtained the large and moderate deviations for the bivariate realized co-volatility; the moderate deviation of jump robust realized volatility was considered by \citet{J2010}. Our paper is the first work considering the large deviation principles of the volatility functionals. We start from the equally spaced sampling scheme, and the results are then extended to the unequally spaced observations.

The remaining of the paper is organized as follows. Section 2 is the setup. In Section 3, we derive the large and moderate deviation principles, for a fixed $u\in \mathbb{R}^+$. In Section 4, we consider the function level (in terms of $u$) large and moderate deviation principles. The unequally spaced sampling case is studied in Section 5. All technical proofs are postponed to the last section.

\section{Setup}\label{sec2}
We consider a process $\{X_t, t\geq0\}$ on some filtered probability space $(\Omega,\mathcal{F},\{\mathcal F_t\}_{t\geq0}, \mathbb{P})$. It is well-known that under the no-arbitrage assumption the efficient price process must follow a semi-martingale, see \citet{DS1994}. We assume that the logarithm of asset price $\{X_t, t\geq0\}$ follows a one-dimensional Brownian semi-martingale of the form
\begin{equation}\label{model}
 \mathrm{d}X_t=a_t\mathrm{d}t+\sigma_t\mathrm{d}W_t+\int_{\mathbb R}\delta(t-, x)\mu(\mathrm{d}t,\mathrm{d}x),
\end{equation}
where, $W$ is a standard Brownian motion, $\{a_t, t\geq 0\}$ and the spot volatility process $\{\sigma_t, t\geq 0\}$ are c$\grave{\text a}$dl$\grave{\text a}$g processes, $\mu(\mathrm{d}t, \mathrm{d}x)$ is a homogeneous Poisson measure with compensator $\mathrm{d}t\otimes \nu(\mathrm{d}x)$, $\nu$ is the L\'{e}vy measure satisfying $\int_{\mathbb R}(x^2\wedge 1)\nu(\mathrm{d}x)<\infty$,
%$\widetilde{\mu}(\mathrm{d}t,\mathrm{d}x):=\mu(\mathrm{d}t,\mathrm{d}x)-\mathrm{d}t\otimes\nu(\mathrm{d}x)$ is the compensated Poisson random measure,
and $\delta(t,x): {\mathbb R}^+\times {\mathbb R}\rightarrow {\mathbb R}$ is c$\grave{\text a}$dl$\grave{\text a}$g in $t$. All the coefficients processes are locally bounded.

The process given in (\ref{model}) is rather a canonical model in the financial theory. Since $\{\sigma_t, t\geq 0\}$ is c$\grave{\text a}$dl$\grave{\text a}$g, all powers of $\sigma$ are locally integrable with respect to the Lebesgue measure. Moreover, we assume that $\{\sigma_t,t\geq 0\}$ is independent of $W$, and bounded below from zero.
%Moreover, both $\{a_s, s\geq 0\}$ and $\{\sigma_s, s\geq 0\}$ can have jumps, intra-day seasonality and long memory.

We need the following assumption on the L\'{e}vy measure due to technical reason.
%The first is about the L\'{e}vy process $L$, and the second imposes some restriction to the smoothness of coefficient processes.%, the assumptions have been used  in Todorov \& Tauchen (2012a).
\begin{asu} \label{asu}
The L\'{e}vy measure $\nu$ satisfies
\begin{eqnarray*}
\mathbb{E}\Big[\int_{0}^{t}\int_{\mathbb R}(|\delta(s, x)|^p \vee |\delta(s, x)|) \nu(\mathrm{d}x)\mathrm{d}s\Big]<\infty,
\end{eqnarray*}
for every $p\in (\beta,1)$, where $0\leq \beta < 1$ is a constant.
\end{asu}
\begin{remark}Assumption \ref{asu} restricts the jump component of $X$, which is required to be with finite variation. The assumption indicates
\begin{equation*}
\mathbb{E}\Big[\int_0^t\int_{|\delta|\geq1}|\delta(s,x)|\nu(\mathrm{d}x)\mathrm{d}s\Big]<\infty,~~ \mathbb{E}\Big[\int_0^t\int_{|\delta|<1}|\delta(s,x)|^p\nu(\mathrm{d}x)\mathrm{d}s\Big]<\infty,
\end{equation*}
for $p\in(\beta, 1)$, which implies that the large jumps are integrable with respect to $\nu(\mathrm{d}x)$ and the small jumps have a Blumenthal-Getoor index $\beta<1$.
Throughout the paper, let $\mathbf P(\cdot)$ and $\mathbf E[\cdot]$ be the conditional probability and conditional expectation given the path of $\sigma$, respectively.
In the proofs, we will assume that the volatility process $\sigma$ is deterministic. If $\sigma$ is random and independent of $W$, the proofs can be easily replicated by conditioning on the path of volatility process, this technique has been used in \citet{KO2011}.
\end{remark}

Recall that the quantity of interest is the realized Laplace transform of volatility: $\int_0^t\exp\{-u\sigma^2_s\}\rm{d}s$. We observe the process $X$ in time interval $[0, t]$ at time points $\{t_i, i=0, 1,\cdots\},$ where $t_i=i\Delta_n$, $\Delta_n$ is the length of the high-frequency interval and $t$ is the time span. Let $t=1$ for simplicity. For $u>0$, the empirical realized Laplace transform of volatility (ERLTV) is defined by \citet{TT2012a} as
\begin{equation}\label{realized Laplace transform}
V_n(u)=\frac{1}{n}\sum_{i=1}^{n}\cos(\sqrt{2nu}\Delta_i^nX),\quad \Delta_i^nX=X_{i/n}-X_{(i-1)/n},
\end{equation}
where $n=1/\Delta_n$. Under the Assumptions A and B in \citet{TT2012a}, the following convergence holds:
\begin{equation}\label{CLT}
\sqrt{n}\Big(V_n(u)-\int_0^1e^{-u\sigma_s^2}\mathrm{d}s\Big)\xrightarrow{\mathcal L-s}\Psi(u),
\end{equation}
where, $\Psi(u)$ is a mixed normal distribution, with $\mathcal F$-conditional variance
$$\int_0^1\frac{e^{-4u\sigma_s^2}-2e^{-2u\sigma_s^2}+1}{2}\mathrm{d}s.$$

The central limit theorem above indicates that an asymptotic approximation to the estimation error between the $V_n(u)$ and RLTV, is given by a (conditional) normal distribution. However, if we would study the tail behaviour of probability distribution of $V_n(u)$, the central limit theorem would not be informative. Therefore, we derive the large deviation principles of $V_n(u)$. More specifically, we are interested in the estimation of
$$\mathbf P\left(m_nV(u)\in A\right),$$
where, $A$ is a given domain of deviation. This includes three types deviations (with different limiting behaviours):
\begin{itemize}
\item if $m_n=\sqrt{n}$, it is just the central limit theorem;
\item if $m_n=1$, it is the large deviation principle;
\item if $m_n\rightarrow+\infty$ and $m_n=o(\sqrt{n})$, it is the moderate deviation principle (after centralization).
\end{itemize}

We consider the large deviation principle in the sense of \citet{DZ1998}. A sequence of $\{Z_n\}_{n\in\mathbb N}$ of random variables satisfies the large deviation principle (LDP) with speed $s_n\downarrow0$ and good rate function $I(\cdot)$, if
\begin{itemize}
  \item $I(\cdot)$ is a good rate function, i.e., $I(\cdot)$ is a lower semicontinuous nonnegative function and the level sets $I^{-1}([0,c])$ are compact for all $c>0$.
  \item (Upper bound) for any closed set $F\in\mathcal B(\mathbb R)$,
  \begin{equation}\label{upper bound}
  \limsup_{n\rightarrow\infty}s_n\log\mathbf P(Z_n\in F)\leq-\inf_{x\in F}I(x),
  \end{equation}
  \item (Lower bound) for any open set $G\in\mathcal B(\mathbb R)$,
  \begin{equation}\label{lower bound}
  \limsup_{n\rightarrow\infty}s_n\log\mathbf P(Z_n\in G)\geq-\inf_{x\in G}I(x).
  \end{equation}
\end{itemize}
We say that $\{Z_n\}_{n\in\mathbb N}$ satisfies the weak large deviation principle with speed $s_n\downarrow0$ and good rate function $I(\cdot)$, if $I(\cdot)$ is a good rate function, the upper bound \eqref{upper bound} holds for all compact sets and the lower bound \eqref{lower bound} holds for all open sets.

 \section{Large and moderate deviations of $V_n(u)$: fixed $u$}

Our first result is about the large deviation principle of $V_n(u)$.
\begin{theorem}\label{LDP}
Assumption \ref{asu} holds. If the volatility $\sigma_s$ is uniformly continuous, then $\mathbf P(V_n(u)\in \cdot)$ satisfies the LDP on $\mathbb R$ with the speed $s_n=1/n$ and the good rate function $I(x,u)$ given by the Fenchel-Legendre transform of $\Lambda(\lambda,u)$, i.e.,
\begin{equation}\label{rate function}
I(x,u)=\Lambda^*(x,u):=\sup_{\lambda\in\mathbb R}\{\lambda x-\Lambda(\lambda,u)\},
\end{equation}
where \begin{equation}\label{limiting logarithmic generating function}
\Lambda(\lambda,u)=\int_0^1\log\Big\{\int_{\mathbb R}\frac{1}{\sqrt{4\pi u\sigma_s^2}}\exp\big(-\frac{y^2}{4u\sigma_s^2}+\lambda\cos y\big)\mathrm{d}y\Big\}\mathrm{d}s.
\end{equation}
\end{theorem}
\begin{rmk}
The good rate function $I(x,u)$ is equal to infinity if $x$ lies beyond an interval containing zero. Figure \ref{fig1} gives a graphical illustration for $I(x,u)$. The curves plot $\lambda x-\Lambda(\lambda,u)$ against $\lambda$, we choose $x$ from -1 to 1 with step 0.05. Hence, totally we obtain 21 curves, where we take $u=1$ and $\sigma_s\equiv1$, the integral is approximated by Riemann sum with fine enough partition. From the left panel, we see that when $x$ is in the interval $[-0.8, 0.8]$, the supremum of $\lambda x-\Lambda(\lambda,u)$ is achieved for finite $\lambda$, whereas for
$x>0.8$, it is $+\infty$ with $\lambda=+\infty$ and $x<-0.8$, it is $+\infty$ with $\lambda=-\infty$.
%the $x$ beyond the interval $[-0.8, 0.8]$,
This fact is further confirmed by the plot in right panel.
\begin{figure}[!ht]
\centering
\includegraphics[width=8cm,height=6cm]{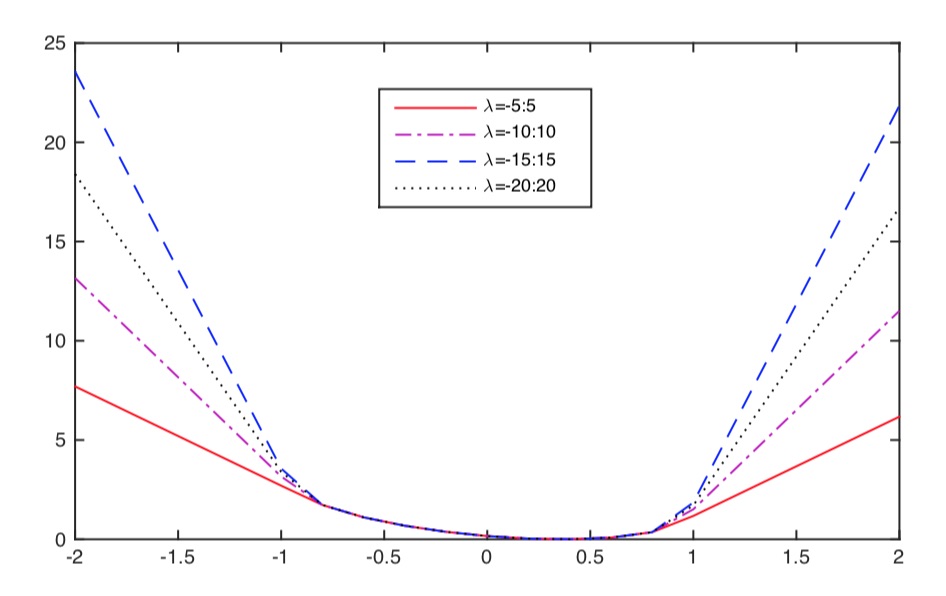}\includegraphics[width=8cm,height=6cm]{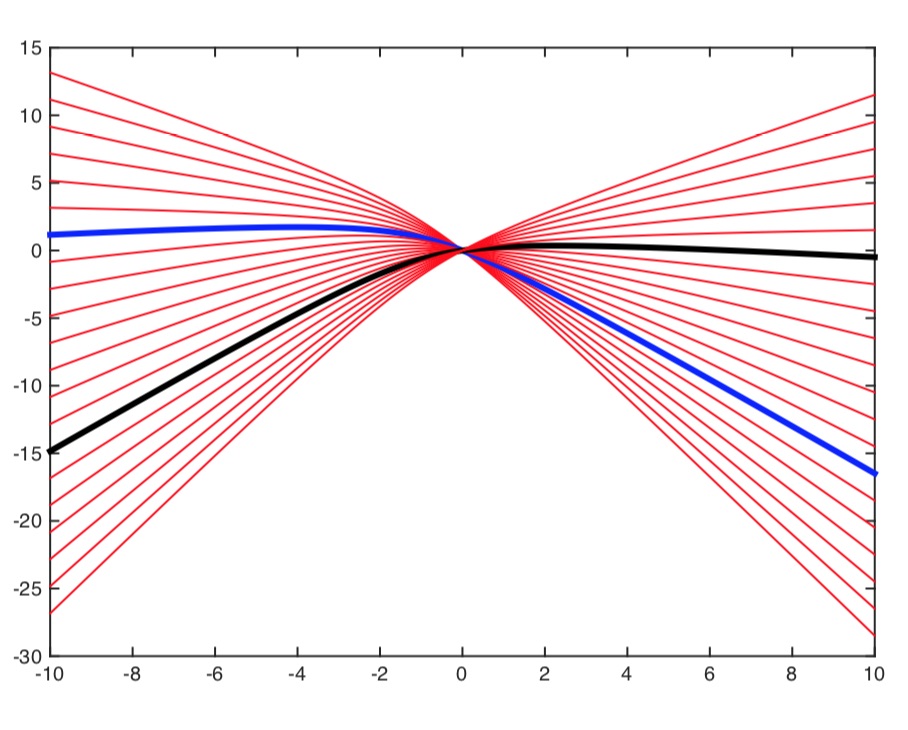}
\caption{Left panel: the $x$-axis is $\lambda$, $y$-axis is $\lambda x-\Lambda(\lambda,u)$. Two bold curves specify the bounds so that $I(x,u)$ be finite, which represents the curve $\lambda x-\Lambda(\lambda,u)$ for $x=-0.8$ and $x=0.8$, respectively. Right panel: the $x$-axis is $x$, $y$-axis is $I(x,u)$, the curves represent the plots of $I(x,u)$ under different range of $\lambda$.}\label{fig1}
\end{figure}
\end{rmk}
\begin{rmk}
From the definition of $V_n(u)$, we know that $|V_n(u)|\leq1$ for all $u\in[0,1]$. Therefore, for $|x|>1$, $\mathbf P(|V_n(u)|>|x|)=0$. Thus, it is necessary that $\inf_{|x|>1}I(x,u)=+\infty$ for the lower bound of LDP to hold. The next proposition verifies this.
\end{rmk}
\begin{proposition}\label{property of rate function}
The limiting logarithmic moment generating function $\Lambda(\lambda,u)$ in Theorem \ref{LDP} satisfies
$$|\Lambda(\lambda,u)|\leq|\lambda|,\quad |\Lambda'(\lambda,u)|\leq1,\quad \Lambda''(\lambda,u)>0,\quad \lambda\in\mathbb{R}.$$
Therefore, the Fenchel-Legedre transform $I(x,u)$ of $\Lambda(\lambda,u)$ is uniquely achieved at some point $\lambda_x$ for all $x\in\mathbb R$. Moreover, for $|x|>1$, we have $I(x,u)=+\infty$.
\end{proposition}

The result of next theorem gives the moderate principle of $V_n(u)$ with fixed $u$.
%\subsection{Moderate Deviation Principle at fixed time}
Let $m_n=o(n^{1/2})$ be a sequence of positive numbers such that $m_n\rightarrow+\infty$ and $m_n/\sqrt{n}\rightarrow0$. We discuss now the moderate deviation of $V_n(u)$.

\begin{theorem}\label{MDP}
Assumption \ref{asu} holds. For $m_n=o(n^{1/2})$, if volatility $\sigma_s$ is $1/2$-H\"{o}lder continuous, then
$$\mathbf P\Big(m_n(V_n(u)-\int_0^1e^{-u\sigma_s^2}\mathrm{d}s)\in \cdot\Big)$$
satisfies the MDP on $\mathbb R$ with the speed $s_n=m_n^2/n$ and the good rate function $$\Lambda^*(x,u)=\frac{x^2}{1+\int_0^1e^{-4u\sigma_s^2}\mathrm{d}s-2\int_0^1e^{-2u\sigma_s^2}\mathrm{d}s}.$$
 \end{theorem}

\section{Large and moderate deviations of $V_n(\cdot)$}
In the previous section, we establish the large and moderate principles for $V_n(u)$ when $u$ is fixed. We next extend the large deviation principle to the function-level, i.e., $V_n(u)$ is in the space of $\mathbb R^{[0,1]}$ of functions on $[0,1]$ indexed by $u$. We define $$\mathcal C([0,1];\mathbb R):=\{\phi:[0,1]\rightarrow\mathbb R|\ \phi(u)\text{ is continuous}\},\quad\|\phi\|:=\sup_{0\leq u\leq1}|\phi(u)|.$$
\begin{theorem}\label{LDP of process}
Assumption \ref{asu} holds. For $m_n=o(n^{1/2})$, if the volatility $\sigma_s$ is uniformly continuous, then $\mathbf P(V_n(\cdot)\in \cdot)$ satisfies the LDP on $(\mathcal C([0,1];\mathbb R),\|\cdot\|)$ with the speed $s_n=1/n$ and some rate function $I^\infty(\cdot)$. Moreover,
$$I^\infty(\phi)\geq\int_0^1I(\phi(u),u)\mathrm{d}u,$$
where $I(\cdot,\cdot)$ is defined in \eqref{rate function}.
\end{theorem}

Next theorem gives the moderate deviation principle of $V_n(u)$.
For $0\leq u\leq 1$, define $F(u):=\int_0^1e^{-u\sigma_s^2}\mathrm{d}s.$
\begin{theorem}\label{MDP of process}
Assumption \ref{asu} holds. For $m_n=o(n^{1/2})$, if the volatility $\sigma_s$ is $1/2$-H\"{o}lder continuous, then
$$\mathbf P\left(m_n\left(V_n(\cdot)-F(\cdot)\right)\in \cdot\right)$$
 satisfies the weak LDP on $\mathcal C([0,1];\mathbb R)$ with respect to the point-wise convergence topology with the speed $s_n=m_n^2/n$ and some rate function $I^\infty(\cdot)$.
\end{theorem}

\section{Extension: unequally spaced sampling}
We now extend the results in previous sections to the unequally sampled observations. Instead of $\{i\Delta_n, i=1\cdots, n\}$, we now observe the process at $\{t_{n,i}, i=1,\cdots, N_n\}$. We adopt the setting of \citet{RSW2017} for the sampling scheme.
\begin{asu} \label{asu:time}
We suppose that the sampling scheme is a double array $(t_{n, i})_{0\leq i\leq N_n; n\geq 1}$ such that
\begin{itemize}
\item $0=t_{n,0}<t_{n,1}<\cdots<t_{n, N_n}\leq 1;$
\item for all $n$, $N_n\leq n$;
\item $\sqrt{n}\max|t_{n,i}-t_{n,i-1}|=o(1)$, as $n\to\infty$;
\item $n\sum_{i=1}^n(t_{n,i}-t_{n,i-1})^2=O(1)$, as $n\to\infty$;
%\item $(t_{n, i})_{0\leq i\leq N_n; n\geq 1}$  is $\{{\cal F}_t\}_{0\leq t\leq 1}$-predictable; that is $t_{n,i}$ is ${\cal F}_{t_{n, i-1}}$ measurable, for all $i, n$.
\end{itemize}
\end{asu}
Define
$$
T_n(s)=\frac{1}{n}\sum_{i=1}^{N_n}\frac{1}{t_{n,i}-t_{n,i-1}}[\min(s,t_{n,i})-\min(s, t_{n,i-1})].
$$
\begin{asu} \label{asu:time2}
The function $T_n$ defined above converges to a function $T$ on $[0,1]$, as $n\to\infty$. Moreover, $T$ admits a strictly positive and continuously bounded derivative $T'$ such that $T_n'\to T'$.
\end{asu}
We now redefine the RLTV as
\begin{equation}\label{NRLT}
V_n(u)=\sum_{i=1}^{N_n}\cos(\sqrt{\frac{2u}{t_{n,i}-t_{n,i-1}}}\Delta_i^nX)(t_{n,i}-t_{n,i-1}), \quad \Delta_i^nX=X_{t_{n,i}}-X_{t_{n,i-1}},
\end{equation}
 According to the Lemma 1 of \citet{RSW2017}, the following convergence holds:
 \begin{equation}\label{CLT}
\sqrt{n}\Big(V_n(u)-\int_0^1e^{-u\sigma_s^2}\mathrm{d}s\Big)\xrightarrow{\mathcal L-s}\Psi(u),
\end{equation}
where, $\Psi(u)$ is a mixed normal distribution, with $\mathcal F$-conditional variance
$$\int_0^1\frac{e^{-4u\sigma_s^2}-2e^{-2u\sigma_s^2}+1}{2T'(s)}ds.$$
We have the following result for the large deviation in the setting of unequally spaced sampling, for the case of fixed $u\in \mathbb{R}$.
\begin{theorem}\label{LDP2}
Assume that Assumptions \ref{asu}-\ref{asu:time2}
 hold.
 \begin{itemize}
 \item[1).]	
If the volatility $\sigma_s$ is uniformly continuous, then $\mathbf P(V_n(u)\in \cdot)$ satisfies the LDP on $\mathbb R$ with the speed $s_n=1/n$ and the good rate function $I(x,u)$ given by the Fenchel-Legendre transform of $\Lambda(\lambda,u)$, i.e.,
\begin{equation}\label{rate function1}
I(x,u)=\Lambda^*(x,u):=\sup_{\lambda\in\mathbb R}\{\lambda x-\Lambda(\lambda,u)\},
\end{equation}
where \begin{equation*}\label{limiting logarithmic generating function1}
\Lambda(\lambda,u)=\int_0^1\log\Big\{\int_{\mathbb R}\frac{1}{\sqrt{4\pi u\sigma_s^2}}\exp\big(-\frac{y^2}{4u\sigma_s^2}+\frac{\lambda\cos y}{T'(s)}\big)\mathrm{d}y\Big\}T'(s)\mathrm{d}s.
\end{equation*}
\item[2).] For $m_n=o(n^{1/4})$, if the volatility $\sigma_s$ is Lipschitz continuous, then
$$\mathbf P\Big(m_n(V_n(u)-\int_0^1e^{-u\sigma_s^2}\mathrm{d}s)\in \cdot\Big)$$
satisfies the MDP on $\mathbb R$ with the speed $s_n=m_n^2/n$ and the good rate function $$\Lambda^*(x,u)=\frac{x^2}{\int_0^1\frac{1+e^{-4u\sigma_s^2}-e^{-2u\sigma_s^2}}{T'(s)}\mathrm{d}s}.$$
\end{itemize}
\end{theorem}

Finally, the large and moderate deviation results of the functional level are given in the theorem below.

\begin{theorem}\label{LDP3}
Assume that Assumptions \ref{asu}-\ref{asu:time2} holds.
\begin{itemize}
\item[1).] If the volatility $\sigma_s$ is uniformly continuous, then $\mathbf P(V_n(\cdot)\in \cdot)$ satisfies the LDP on $(\mathcal C([0,1];\mathbb R),\|\cdot\|)$ with the speed $s_n=1/n$ and some rate function $I^\infty(\cdot)$. Moreover,
$$I^\infty(\phi)\geq\int_0^1I(\phi(u),u)\mathrm{d}u,$$
where $I(\cdot,\cdot)$ is defined in \eqref{rate function1}.
\item[2).] For $m_n=o(n^{1/4})$, if the volatility $\sigma_s$ is Lipschitz continuous, then
$$\mathbf P\left(m_n\left(V_n(\cdot)-F(\cdot)\right)\in \cdot\right)$$
 satisfies the weak LDP on $\mathcal C([0,1];\mathbb R)$ with respect to the point-wise convergence topology with the speed $s_n=m_n^2/n$ and some rate function $I^\infty(\cdot)$.
 \end{itemize}
\end{theorem}

\begin{rmk}
Similar to the proof of Theorem \ref{MDP}, in order to get MDP for unequally spaced sampling, we need to show that
\begin{equation}\nonumber
\begin{aligned}
 \lim_{n\rightarrow\infty}\lambda m_n\sum_{i=1}^{N_n}\int_{t_{n,i-1}}^{t_{n,i}}\Big(e^{-u\sigma_{t_{n,i-1}}^2}-e^{-u\sigma_{s}^2}\Big)\mathrm{d}s
\leq\lim_{n\rightarrow\infty}C\lambda m_n\sum_{i=1}^{N_n}\int_{t_{n,i-1}}^{t_{n,i}}|\sigma_{t_{n,i-1}}-\sigma_{s}|\mathrm{d}s
=0,
\end{aligned}
\end{equation}
and
\begin{equation}\nonumber
\begin{aligned}
&\lim_{n\rightarrow\infty}\frac{m_n^2}{n}\sum_{i=1}^{N_n}\mathbf E \Big[\exp\Big\{\frac{n\lambda\Delta_i^nt}{m_n}\cos(\sqrt{2u/\Delta_i^nt}\Delta_i^nX)\Big\}-
\exp\Big\{\frac{n\lambda\Delta_i^nt}{m_n}\cos(\sqrt{2u/\Delta_i^nt}\sigma_{t_{n,i-1}}\Delta_i^nW)\Big\}\Big]\\
&\leq\lim_{n\rightarrow\infty}\frac{m_n^2}{n}\sum_{i=1}^{N_n}\frac{C}{m_n}|\Delta_i^nt|^{1/2}
\leq\lim_{n\rightarrow\infty}C\frac{m_n}{n}N_n\max_i|\Delta_i^nt|^{1/2}
=0,
\end{aligned}
\end{equation}
where $\Delta_i^n t =t_{n,i}-t_{n,i-1}$. Therefore, the conditions for MDP in Theorem \ref{LDP2} and Theorem \ref{LDP3} are stronger than those in Theorem \ref{MDP}.
\end{rmk}

\section{Proofs}
In the proofs, $c$ denotes a constant independent of $n$ whose value may change from line to line.

\textbf{Proof of Theorem \ref{LDP}:}
For $n\in\mathbb N$, $u\in[0,1]$ and $\lambda\in\mathbb R$, we define the logarithmic moment generating function of the realized volatility statistic as $$\Lambda_n(\lambda,u):=\log\mathbf E\Big[\exp\Big\{\frac \lambda n\sum_{i=1}^n\cos(\sqrt{2nu}\Delta_i^nX)\Big\}\Big].$$
It is easy to see that $|\Lambda(\lambda,u)|\leq|\lambda|$ and $\Lambda(\lambda,u)$ is differentiable with respect to $\lambda$ on $\mathbb R$. Therefore, according to the G\"{a}rtner-Ellis's Theorem in \citet{DZ1998}, it suffices to show that for each $\lambda\in\mathbb R$,
$$\lim_{n\rightarrow\infty}\frac 1 n \Lambda_n(\lambda n,u)=\Lambda(\lambda,u).$$
 It follows from the independent increments of $X$ that
\begin{equation}\label{501}
\begin{aligned}
\frac 1 n\Lambda_n(\lambda n,u)=&\frac 1 n \log\mathbf E\Big[\exp\Big\{\lambda\sum_{i=1}^n\cos(\sqrt{2nu}\Delta_i^nX)\Big\}\Big]\\
=&\frac 1 n\sum_{i=1}^n\log\mathbf E\Big[\exp\Big\{\lambda\cos(\sqrt{2nu}\Delta_i^nX)\Big\}\Big]=\frac 1 n\sum_{i=1}^n\log\mathbf E\big[A_{n1}+A_{n2}\big],\\
\end{aligned}
\end{equation}
where $A_{n1}:=\exp\Big\{\lambda\cos(\sqrt{2 n u }\sigma_{t_{i-1}}\Delta_i^nW)\Big\},$
and $$A_{n2}:=\exp\Big\{\lambda\cos(\sqrt{2nu}\Delta_i^nX)\Big\}-\exp\Big\{\lambda\cos(\sqrt{2nu}\sigma_{t_{i-1}}\Delta_i^nW)\Big\}.$$
First, note that
\begin{equation}\label{A_{n1}}
\begin{aligned}
&\mathbf E[A_{n1}]
=&
\int_{\mathbb R}\frac{1}{\sqrt{4\pi u\sigma_{t_{i-1}}^2}}\exp\Big\{-\frac{y^2}{4u\sigma_{t_{i-1}}^2}+\lambda\cos y\Big\}\mathrm{d}y
:=&\phi_{i,n}(\lambda,u).\\
\end{aligned}
\end{equation}

Next,
\begin{eqnarray*}
\mathbf E|A_{n2}|
&\leq&\mathbf E\Big[|\lambda|e^{|\lambda|}\big|\cos(\sqrt{2 n u}\Delta_i^nX)-\cos(\sqrt{2 n u}\sigma_{(i-1)/n}\Delta_i^nW)\big|\Big]\\
&\leq&|\lambda|e^{|\lambda|}\mathbf E\big[|B_{n1}|+|B_{n2}|\big],
\end{eqnarray*}
where
\begin{equation}\nonumber
\begin{aligned}
B_{n1}
=&-2\sin\Big(\sqrt{2nu}\big(\int_{t_{i-1}}^{t_i}a_s\mathrm{d}s+\int_{t_{i-1}}^{t_i}\sigma_s\mathrm{d}W_s
+\frac{1}{2}\int_{t_{i-1}}^{t_i}\int_{\mathbb R}\delta(s-, x)\mu(\mathrm{d}s,\mathrm{d}x)\big)\Big)\\
&\sin\Big(\frac{\sqrt{2nu}}{2}\int_{t_{i-1}}^{t_i}\int_{\mathbb R}\delta(s-, x)\mu(\mathrm{d}s,\mathrm{d}x)\Big),
\end{aligned}
\end{equation}
and
\begin{eqnarray*}
B_{n2}=\cos\Big(\sqrt{2nu}\big(\int_{t_{i-1}}^{t_i}a_s\mathrm{d}s+\int_{t_{i-1}}^{t_i}\sigma_s\mathrm{d}W_s\big)\Big)
-\cos\Big(\sqrt{2nu}\sigma_{t_{i-1}}\Delta_i^nW\Big).
\end{eqnarray*}
Next, we study $B_{n1}$ and $B_{n2}$, respectively.
%From Assumption \ref{asu}(2), we have $\Delta_i^nL=^d(1/n)^{1/\alpha}L_1$.
For $\beta<1$, choose $\epsilon>0$ such that $\beta+\epsilon<1$. Therefore, it follows from $|\sin(x)|\leq|\sin(x)|^{\beta+\epsilon}\leq|x|^{\beta+\epsilon}$ and $C_r$-inequality that
\begin{eqnarray*}
&&\mathbf E\Big|\sin\Big(\frac{\sqrt{2nu}}{2}\int_{t_{i-1}}^{t_i}\int_{\mathbb R}\delta(s-, x)\mu(\mathrm{d}s,\mathrm{d}x)\Big)\Big|\\
&\leq&\Big(\frac{\sqrt{2nu}}{2}\Big)^{\beta+\epsilon}\mathbf E\Big|\int_{t_{i-1}}^{t_i}\int_{\mathbb R}\delta(s-, x)\mu(\mathrm{d}s,\mathrm{d}x)\Big|^{\beta+\epsilon}\\
&\leq&\Big(\frac{\sqrt{2nu}}{2}\Big)^{\beta+\epsilon}\mathbf E\int_{t_{i-1}}^{t_i}\int_{\mathbb R}\big|\delta(s-, x)\big|^{\beta+\epsilon}\nu(\mathrm{d}x)\mathrm{d}s\\
&\leq&\frac{c}{n^{1-\beta/2-\epsilon/2}}.
\end{eqnarray*}
Therefore,
\begin{eqnarray}\label{B_{n1}}
\mathbf E|B_{n1}|=o\big(\frac{1}{n^{1/2}}\big).
\end{eqnarray}
It follows the mean value theorem that
\begin{eqnarray*}
B_{n2}
&=&\cos\left(\sqrt{2 n u}\int_{t_{i-1}}^{t_i}a_s\mathrm{d}s+\sqrt{2nu}
\int_{t_{i-1}}^{t_i}\sigma_s\mathrm{d}W_s\right)-\cos\Big(\sqrt{2u}n^{-\frac{1}{2}}a_{t_{i-1}}+
\sqrt{2nu}\int_{t_{i-1}}^{t_i}\sigma_s\mathrm{d}W_s\Big)\\
&&+\cos\Big(\sqrt{2u}n^{-\frac{1}{2}}a_{t_{i-1}}+
\sqrt{2nu}\int_{t_{i-1}}^{t_i}\sigma_s\mathrm{d}W_s\Big)-\cos\Big(\sqrt{2 n u}\sigma_{t_{i-1}}
\Delta_i^nW\Big)\\
&=&-\Big(\sqrt{2nu}\int_{t_{i-1}}^{t_i}a_s\mathrm{d}s-
\sqrt{2u}n^{-\frac{1}{2}}a_{t_{i-1}}\Big)\sin\big(\widetilde{x}_1\big)\\
&&-\sqrt{2nu}
\Big(n^{-1}a_{t_{i-1}}+\int_{t_{i-1}}^{t_i}(\sigma_s-\sigma_{t_{i-1}})\mathrm{d}W_s\Big)\sin\big(\widetilde{x}_2\big)\\
&=&-\Big(\sqrt{2nu}\int_{t_{i-1}}^{t_i}a_s\mathrm{d}s-
\sqrt{2u}n^{-\frac{1}{2}}a_{t_{i-1}}\Big)\sin\big(\widetilde{x}_1\big)-\sqrt{2u}n^{-\frac{1}{2}}a_{t_{i-1}}\sin\Big(\sqrt{2 n u}\sigma_{t_{i-1}}\Delta_i^nW\Big)
\\
&&-\sqrt{2nu}\sin\Big(\sqrt{2 n u}\sigma_{t_{i-1}}\Delta_i^nW\Big)
\int_{t_{i-1}}^{t_i}(\sigma_s-\sigma_{t_{i-1}})\mathrm{d}W_s\\
&&+\sqrt{2nu}\Big(\frac{1}{n}a_{t_{i-1}}+\int_{t_{i-1}}^{t_i}(\sigma_s-\sigma_{t_{i-1}})\mathrm{d}W_s\Big)\Big(\sqrt{2nu}
\sigma_{t_{i-1}}\Delta_i^nW-\widetilde{x}_2\Big)\cos\big(\widetilde{x}_3\big)\\
&:=&\xi_1+\xi_2+\xi_3+\xi_4,
\end{eqnarray*}
where
\begin{itemize}
  \item $\widetilde{x}_1$ is between $\sqrt{2 n u}\int_{t_{i-1}}^{t_i}a_s\mathrm{d}s+\sqrt{2nu}\int_{t_{i-1}}^{t_i}\sigma_s\mathrm{d}W_s$ and $\sqrt{2u}n^{-\frac{1}{2}}a_{t_{i-1}}+\sqrt{2nu}\int_{t_{i-1}}^{t_i}\sigma_s\mathrm{d}W_s$;
  \item $\widetilde{x}_2$ is between $\sqrt{2nu}\int_{t_{i-1}}^{t_i}a_s\mathrm{d}s+\sqrt{2nu}\int_{t_{i-1}}^{t_i}\sigma_s\mathrm{d}W_s$
and $\sqrt{2nu}\sigma_{t_{i-1}}\Delta_i^nW$;
  \item $\widetilde{x}_3$ is between $\widetilde{x}_2$ and $\sqrt{2nu}\sigma_{t_{i-1}}\Delta_i^nW$.
\end{itemize}
It is easy to check that
\begin{equation}\nonumber
\begin{aligned}
\mathbf E|\xi_1|\leq cn^{1/2}\mathbf E\int_{t_{i-1}}^{t_i}|a_s-a_{t_{i-1}}|\mathrm{d}s\leq cn^{-1/2},
\end{aligned}
\end{equation}
and
\begin{equation}
\mathbf E|\xi_2|\leq cn^{1/2}\mathbf E\int_{t_{i-1}}^{t_i}|a_{t_{i-1}}|\mathrm{d}s\leq cn^{-1/2}.
\end{equation}
It the volatility $\sigma_s$ is uniformly continuous, then there exist a sequence $\varepsilon_n\rightarrow0$ as $n\rightarrow\infty$ such that
\begin{equation}\nonumber
\begin{aligned}
\mathbf E|\xi_3|\leq cn^{1/2}\Big(\mathbf E\int_{t_{i-1}}^{t_i}|\sigma_s-\sigma_{t_{i-1}}|^2\mathrm{d}s\Big)^{1/2}\leq c\varepsilon_n,
\end{aligned}
\end{equation}
and
\begin{equation}\nonumber
\begin{aligned}
\mathbf E|\xi_4|\leq &cn^{1/2}\Big(\mathbf E\Big(\frac{1}{n}+\int_{t_{i-1}}^{t_i}(\sigma_s-\sigma_{t_{i-1}})\mathrm{d}W_s\Big)^2\Big)^{1/2}
\\&\Big(\mathbf E\Big(\sqrt{2nu}\int_{t_{i-1}}^{t_i}(\sigma_{t_{i-1}}-\sigma_s)\mathrm{d}W_s-\sqrt{2nu}\int_{t_{i-1}}^{t_i}a_s\mathrm{d}s\Big)^2\Big)^{1/2}\\
\leq&cn\Big(\frac{1}{n^2}+\mathbf E\int_{t_{i-1}}^{t_i}(\sigma_s-\sigma_{t_{i-1}})^2\mathrm{d}s\Big)\\
\leq&c\varepsilon_n^2.
\end{aligned}
\end{equation}
Therefore,
\begin{eqnarray}\label{B_{n2}}
\mathbf E|B_{n2}|\leq c\varepsilon_n.
\end{eqnarray}
Combing \eqref{B_{n1}} and \eqref{B_{n2}}, we have
\begin{eqnarray}\label{A_{n2}}
\mathbf E|A_{n2}|\leq c\varepsilon_n+c\frac{1}{n^{1/2}}.
\end{eqnarray}
Therefore, for large enough $n$, it follows from \eqref{501}, \eqref{A_{n1}} and \eqref{A_{n2}} that
\begin{equation}\nonumber
\begin{aligned}
\frac 1 n \sum_{i=1}^n\Lambda_n(\lambda n,u)
\leq&\frac 1 n\sum_{i=1}^n\log\Big(\phi_{i,n}(\lambda,u)+\frac{c}{n^{1/2}}+c\varepsilon_n\Big)\\
\leq&\frac 1 n\sum_{i=1}^n\log(\phi_{i,n}(\lambda,u))+\frac{ce^{|\lambda|}}{n^{1/2}}+c\varepsilon_ne^{|\lambda|},
\end{aligned}
\end{equation}
and
\begin{equation}\nonumber
\begin{aligned}
\frac 1 n \sum_{i=1}^n\Lambda_n(\lambda n,u)
\geq&\frac 1 n\sum_{i=1}^n\log\Big(\phi_{i,n}(\lambda,u)-c\frac{1}{n^{1/2}}-c\varepsilon_n\Big)\\
\geq&\frac 1 n\sum_{i=1}^n\log(\phi_{i,n}(\lambda,u))
-\frac{2ce^{-|\lambda|}}{n^{1/2}}-2ce^{-|\lambda|}.
\end{aligned}
\end{equation}

Thus, by Dominated Convergence Theorem, we have
\begin{equation}\nonumber
\begin{aligned}
&\limsup_{n\rightarrow\infty}\frac 1 n \sum_{i=1}^n\Lambda_n(\lambda n,u)\\
\leq&\limsup_{n\rightarrow\infty}\frac 1 n\sum_{i=1}^n\log(\phi_{i,n}(\lambda,u))\\
=&\limsup_{n\rightarrow\infty}\int_0^1\sum_{i=1}^n1_{[(i-1)/n,i/n)}(s)\log\int_{\mathbb R}\frac{1}{\sqrt{4\pi u\sigma_{(i-1)/n}^2}}\exp\Big\{-\frac{y^2}{4u\sigma_{(i-1)/n}^2}+\lambda\cos y\Big\}\mathrm{d}y\mathrm{d}s\\
=&\int_0^1\log\int_{\mathbb R}\frac{1}{\sqrt{4\pi u\sigma_s^2}}\exp\Big\{-\frac{y^2}{4u\sigma_s^2}+\lambda\cos y\Big\}\mathrm{d}y\mathrm{d}s,\\
\end{aligned}
\end{equation}
and
\begin{eqnarray*}
\liminf_{n\rightarrow\infty}\frac 1 n \sum_{i=1}^n\Lambda_n(\lambda n,u)
&\geq&\liminf_{n\rightarrow\infty}\frac 1 n\sum_{i=1}^n\log(\phi_{i,n}(\lambda,u))\\
&=&\int_0^1\log\int_{\mathbb R}\frac{1}{\sqrt{4\pi u\sigma_s^2}}\exp\Big\{-\frac{y^2}{4u\sigma_s^2}+\lambda\cos y\Big\}\mathrm{d}y\mathrm{d}s.\\
\end{eqnarray*}
Therefore, for $\lambda\in\mathbb R$,
$$\hspace{5cm} \lim_{n\rightarrow\infty}\frac 1 n \sum_{i=1}^n\Lambda_n(\lambda n,u)=\Lambda(\lambda,u).\hspace{5cm} \Box$$

\

\textbf{Proof of Proposition \ref{property of rate function}:}
From the definition of $\Lambda(\lambda,u)$, we have
$$\Lambda'(\lambda,u)=\int_0^1\frac{\int_\mathbb R\frac{1}{\sqrt{4\pi u\sigma_s^2}}\exp\Big\{-\frac{y^2}{4u\sigma_s^2}+\lambda\cos y\Big\}\cos y\mathrm{d}y}{\int_\mathbb R\frac{1}{\sqrt{4\pi u\sigma_s^2}}\exp\Big\{-\frac{y^2}{4u\sigma_s^2}+\lambda\cos y\Big\}\mathrm{d}y}\mathrm{d}s,$$
and
\begin{equation}\nonumber
\begin{aligned}
\Lambda''(\lambda,u)=&\int_0^1\frac{\int_\mathbb R\frac{1}{\sqrt{4\pi u\sigma_s^2}}\exp\Big\{-\frac{y^2}{4u\sigma_s^2}+\lambda\cos y\Big\}\cos^2 y\mathrm{d}y\int_\mathbb R\frac{1}{\sqrt{4\pi u\sigma_s^2}}\exp\Big\{-\frac{y^2}{4u\sigma_s^2}+\lambda\cos y\Big\}\mathrm{d}y}{\Big(\int_\mathbb R\frac{1}{\sqrt{4\pi u\sigma_s^2}}\exp\Big\{-\frac{y^2}{4u\sigma_s^2}+\lambda\cos y\Big\}\mathrm{d}y\Big)^2}\mathrm{d}s\\
&-\int_0^1\frac{\Big(\int_\mathbb R\frac{1}{\sqrt{4\pi u\sigma_s^2}}\exp\Big\{-\frac{y^2}{4u\sigma_s^2}+\lambda\cos y\Big\}\cos y\mathrm{d}y\Big)^2}{\Big(\int_\mathbb R\frac{1}{\sqrt{4\pi u\sigma_s^2}}\exp\Big\{-\frac{y^2}{4u\sigma_s^2}+\lambda\cos y\Big\}\mathrm{d}y\Big)^2}\mathrm{d}s
\end{aligned}
\end{equation}
Therefore, it follows from $|\cos(y)|\leq1$ and the H\"{o}lder inequality that for $\lambda\in\mathbb R$,
$$|\Lambda(\lambda,u)|\leq|\lambda|,\quad |\Lambda'(\lambda,u)|\leq1,\quad \Lambda''(\lambda,u)>0.$$
Since $\Lambda'(\lambda,u)$ is a strictly increasing function, $I(x,u)=\sup_{\lambda\in\mathbb R}\{\lambda x-\Lambda(\lambda,u)\}$ is uniquely achieved at some point. For $x>1$, $x-\Lambda'(\lambda,u)>0$, thus $I(x)$ is achieved as $\lambda\rightarrow+\infty$. Recall $|\Lambda(\lambda,u)|\leq\lambda$, therefore,
$$I(x,u)=+\infty,\qquad x>1.$$
Similarly,
$$\hspace{5cm} I(x,u)=+\infty,\qquad x<-1. \hspace{5cm} \Box$$

\

\textbf{Proof of Theorem \ref{MDP}:}
If $\sigma_s$ is $1/2$-H\"{o}lder continuous, we have
\begin{eqnarray*}
\mathbf E\int_{t_{i-1}}^{t_i}(\sigma_s-\sigma_{t_{i-1}})^2\mathrm{d}s\leq\frac{c}{n^2}.
\end{eqnarray*}
Thus, similar to the computation of $\mathbf E|B_{n2}|$ in Theorem \ref{LDP}, we have
\begin{eqnarray*}
\mathbb E\Big[\Big|\cos(\sqrt{2nu}(\int_{t_{i-1}}^{t_i}a_s\mathrm{d}s+\int_{t_{i-1}}^{t_i}\sigma_s\mathrm{d}W_s))
-\cos(\sqrt{2nu}\sigma_{t_{i-1}}\Delta_i^nW)\Big|\Big]\leq\frac{c}{n^{1/2}}.
\end{eqnarray*}
%For $\beta<1$, choose $\epsilon>0$ such that $\beta+\epsilon<1$. Therefore, it follows from %$|\sin(x)|\leq|\sin(x)|^{\beta+\epsilon}\leq|x|^{\beta+\epsilon}$ and $C_r$-inequality that
%\begin{eqnarray*}
%&&\mathbf E\Big|\sin\Big(\frac{\sqrt{2nu}}{2}\int_{t_{i-1}}^{t_i}\int_{\mathbb R}\delta(s-, x)\mu(\mathrm{d}s,\mathrm{d}x)\Big)\Big|\\
%&\leq&\Big(\frac{\sqrt{2nu}}{2}\Big)^{\beta+\epsilon}\mathbf E\Big|\int_{t_{i-1}}^{t_i}\int_{\mathbb R}\delta(s-, x)\mu(\mathrm{d}s,\mathrm{d}x)\Big|^{\beta+\epsilon}\\
%&\leq&\Big(\frac{\sqrt{2nu}}{2}\Big)^{\beta+\epsilon}\mathbf E\int_{t_{i-1}}^{t_i}\int_{\mathbb R}\big|\delta(s-, x)\big|^{\beta+\epsilon}\mathrm{d}s\nu(\mathrm{d}x)\\
%&\leq&\frac{c}{n^{1-\beta/2-\epsilon/2}}.
%\end{eqnarray*}
It has been shown in Theorem \ref{LDP} that \begin{eqnarray*}
\mathbf E\Big|\sin\Big(\frac{\sqrt{2nu}}{2}\int_{t_{i-1}}^{t_i}\int_{\mathbb R}\delta(s-, x)\mu(\mathrm{d}s,\mathrm{d}x)\Big)\Big|
\leq\frac{c}{n^{1-\beta/2-\epsilon/2}}.
\end{eqnarray*}
Thus,
\begin{eqnarray*}
&&\mathbf E\Big[\exp\Big\{\frac{\lambda}{m_n}\cos(\sqrt{2 n u}\Delta_i^nX)\Big\}-
\exp\Big\{\frac{\lambda}{m_n}\cos(\sqrt{2 n u}\sigma_{t_{i-1}}\Delta_i^nW)\Big\}\Big]\\
&\leq& |\lambda|e^{|\lambda|}m_n^{-1}\mathbf E\Big[\big|\cos(\sqrt{2 n u}\Delta_i^nX)-
\cos(\sqrt{2 n u}\sigma_{t_{i-1}}\Delta_i^nW)\big|\Big]\\
&\leq&|\lambda|e^{|\lambda|}m_n^{-1}\mathbf E\Big[2\Big|\sin\big(\frac{\sqrt{2nu}}{2}\int_{t_{i-1}}^{t_i}\int_{\mathbb R}\delta(s-, x)\mu(\mathrm{d}s,\mathrm{d}x)\big)\Big|\\
&&\qquad+\Big|\cos\Big(\sqrt{2nu}\big(\int_{t_{i-1}}^{t_i}a_s\mathrm{d}s+\int_{t_{i-1}}^{t_i}\sigma_s\mathrm{d}W_s\big)\Big)
-\cos\big(\sqrt{2nu}\sigma_{t_{i-1}}\Delta_i^nW\big)\Big|\Big]\\
&\leq&\frac{c}{m_nn^{1/2}}.
\end{eqnarray*}
From the Taylor expansion
$$e^y=1+y+\frac{y^2}{2}+\frac{y^3}{6}e^{\theta_y},\quad 0<\theta_y<1,$$
we have
\begin{eqnarray*}
&&\frac{m_n^2}{n}\log\mathbf E\Big[\exp\Big\{
\frac{\lambda n}{m_n}\Big(\sum_{i=1}^n\frac 1 n\cos(\sqrt{2nu}\Delta_i^nX)-\int_0^1e^{-u\sigma_s^2}\mathrm{d}s\Big)\Big\}\Big]\\
&=&\frac{m_n^2}{n}\sum_{i=1}^n\log\mathbf E\Big[\exp\Big\{
\frac{\lambda}{m_n}\cos(\sqrt{2nu}\Delta_i^nX)\Big\}\Big]-\lambda m_n\int_0^1e^{-u\sigma_s^2}\mathrm{d}s\\
&\leq&\frac{m_n^2}{n}\sum_{i=1}^n\log\Big[\int_{\mathbb R}\frac{1}{\sqrt{4\pi u\sigma^2_{t_{i-1}}}}\exp\Big\{-\frac{y^2}{4u\sigma^2_{t_{i-1}}}+\frac{\lambda}{m_n}\cos y\Big\}\mathrm{d}y+ \frac{C}{m_nn^{1/2}}\Big]\\
&&-\lambda m_n\int_0^1e^{-u\sigma_s^2}\mathrm{d}s\\
&\leq&\frac{m_n^2}{n}\sum_{i=1}^n\log\big[\int_{\mathbb R}\frac{1}{\sqrt{4\pi u\sigma^2_{t_{i-1}}}}\Big(1+\frac{\lambda}{m_n}\cos y+\frac{\lambda^2\cos^2y}{2m_n^2}+\frac{\lambda^3\cos^3y}{6m_n^3}e^{\frac{\theta_{y,n}\lambda\cos y}{m_n}}\Big)\times\\
&&\exp\Big\{-\frac{y^2}{4u\sigma^2_{t_{i-1}}}\Big\}\mathrm{d}y\big]-\lambda m_n\int_0^1e^{-u\sigma_s^2}\mathrm{d}s+O(\frac{m_ne^{\frac{|\lambda|}{m_n}}}{n^{1/2}}).\\
\end{eqnarray*}
Note that
$$\int_{\mathbb R}\frac{\lambda\cos x}{m_n\sqrt{4\pi u\sigma^2_{t_{i-1}}}}\exp\Big\{-\frac{x^2}{4u\sigma^2_{t_{i-1}}}\Big\}\mathrm{d}x=\frac{\lambda}{m_n}e^{-u\sigma^2_{t_{i-1}}},$$
and
\begin{equation}\nonumber
\begin{aligned}
\int_{\mathbb R}\frac{\lambda^2\cos^2y}{2m_n^2\sqrt{4\pi u\sigma^2_{t_{i-1}}}}\exp\Big\{-\frac{y^2}{4u\sigma^2_{t_{i-1}}}\Big\}\mathrm{d}y=&\frac{\lambda^2}{4m_n^2}
\int_{\mathbb R}\frac{ 1+\cos2y}{\sqrt{4\pi u\sigma^2_{t_{i-1}}}}\exp\Big\{-\frac{y^2}{4u\sigma^2_{t_{i-1}}}\Big\}\mathrm{d}y\\
=&\frac{\lambda^2}{4m_n^2}\Big(1+e^{-4u\sigma^2_{t_{i-1}}}\Big),
\end{aligned}
\end{equation}
and
$$\int_{\mathbb R}\frac{\lambda^3\cos^3y}{6m_n^3\sqrt{4\pi u\sigma^2_{t_{i-1}}}}\exp\Big\{-\frac{y^2}{4u\sigma^2_{t_{i-1}}}+\frac{\theta_{y,n}\lambda\cos y}{m_n}\Big\}\mathrm{d}y\leq\frac{|\lambda|^3e^{|\lambda|}}{6m_n^3},$$
we have
\begin{eqnarray*}
&&\frac{m_n^2}{n}\log\mathbf E\Big[\exp\Big\{
\frac{\lambda n}{m_n}\Big(\sum_{i=1}^n\frac 1 n\cos(\sqrt{2nu}\Delta_i^nX)-\int_0^1e^{-u\sigma_s^2}\mathrm{d}s\Big)\Big\}\Big]\\
&\leq&\frac{m_n^2}{n}\sum_{i=1}^n\log\Big[1+\frac{\lambda}{m_n}e^{-u\sigma^2_{t_{i-1}}}+\frac{\lambda^2}{4m_n^2}\Big(1+e^{-4u\sigma^2_{t_{i-1}}}\Big)
+\frac{|\lambda|^3e^{|\lambda|}}{6m_n^3}\Big]\\
&&+c|\lambda|e^{2|\lambda|}n^{-1/2}m_n-\lambda m_n\int_0^1e^{-u\sigma_s^2}\mathrm{d}s\\
&\leq&\frac{m_n^2}{n}\sum_{i=1}^n\log\Big[1+\frac{\lambda}{m_n}e^{-u\sigma^2_{t_{i-1}}}+\frac{\lambda^2}{4m_n^2}\Big(1+e^{-4u\sigma^2_{t_{i-1}}}\Big)\Big]
\\
&&+\frac{|\lambda|^3e^{|\lambda|}}{6m_n}-\lambda m_n\int_0^1e^{-u\sigma_s^2}\mathrm{d}s+O(\frac{m_ne^{\frac{|\lambda|}{m_n}}}{n^{1/2}}).
\end{eqnarray*}
Thus, it follows from $$\log (1+y)=y-\frac{y^2}{2(1+\eta_yy)^2},\quad 0<\eta_y<1,$$ that
\begin{equation*}
\begin{aligned}
&\frac{m_n^2}{n}\log\mathbf E\Big[\exp\Big\{
\frac{\lambda n}{m_n}\Big(\sum_{i=1}^n\frac 1 n\cos(\sqrt{2nu}\Delta_i^nX)-\int_0^1e^{-u\sigma_s^2}\mathrm{d}s\Big)\Big\}\Big]\\
\leq&\frac{m_n^2}{n}\sum_{i=1}^n\Big[\frac{\lambda}{m_n}e^{-u\sigma^2_{t_{i-1}}}+\frac{\lambda^2}{4m_n^2}\Big(1+e^{-4u\sigma^2_{t_{i-1}}}\Big)\Big]
\\
&-\frac{m_n^2}{2n}\sum_{i=1}^n\frac{\lambda^2}{m_n^2}e^{-2u\sigma^2_{t_{i-1}}}
\Big[1+\eta_{i,n}\Big(\frac{\lambda}{m_n}e^{-u\sigma^2_{t_{i-1}}}+\frac{\lambda^2}{4m_n^2}\Big(1+e^{-4u\sigma^2_{t_{i-1}}}\Big)\Big)\Big]^{-2}\\
&-\frac{m_n^2}{4n}\sum_{i=1}^n\frac{\lambda^3}{m_n^3}e^{-u\sigma^2_{t_{i-1}}}\Big(1+e^{-4u\sigma^2_{t_{i-1}}}\Big)
\Big[1+\eta_{i,n}\Big(\frac{\lambda}{m_n}e^{-u\sigma^2_{t_{i-1}}}+\frac{\lambda^2}{4m_n^2}\Big(1+e^{-4u\sigma^2_{t_{i-1}}}\Big)\Big)\Big]^{-2}\\
&-\frac{m_n^2}{32n}\sum_{i=1}^n\frac{\lambda^4}{m_n^4}\Big(1+e^{-4u\sigma^2_{t_{i-1}}}\Big)^2
\Big[1+\eta_{i,n}\Big(\frac{\lambda}{m_n}e^{-u\sigma^2_{t_{i-1}}}+\frac{\lambda^2}{4m_n^2}\Big(1+e^{-4u\sigma^2_{t_{i-1}}}\Big)\Big)\Big]^{-2}\\
&+\frac{|\lambda|^3e^{|\lambda|}}{6m_n}-\lambda m_n\int_0^1e^{-u\sigma_s^2}\mathrm{d}s+O(\frac{m_ne^{\frac{|\lambda|}{m_n}}}{n^{1/2}})\\
\leq&\lambda m_n\Big(\frac 1 n\sum_{i=1}^ne^{-u\sigma^2_{t_{i-1}}}-\int_0^1e^{-u\sigma_s^2}\mathrm{d}s\Big)
+\frac{\lambda^2}{4n}\sum_{i=1}^n\Big(1+e^{-4u\sigma^2_{t_{i-1}}}\Big)\\
&-\frac{\lambda^2}{2n}\Big(1+\frac{3\lambda^2}{m_n^2}\Big)^{-2}\sum_{i=1}^ne^{-2u\sigma^2_{t_{i-1}}}-\frac{\lambda^3}{4nm_n}
\Big(1+\frac{3\lambda^2}{m_n^2}\Big)^{-2}\sum_{i=1}^ne^{-u\sigma^2_{t_{i-1}}}\Big(1+e^{-4u\sigma^2_{t_{i-1}}}\Big)\\
&-\frac{\lambda^4}{32nm_n^2}\Big(1+\frac{3\lambda^2}{m_n^2}\Big)^{-2}\sum_{i=1}^n\Big(1+e^{-4u\sigma^2_{t_{i-1}}}\Big)^2+\frac{|\lambda|^3e^{|\lambda|}}{6m_n}
+O(\frac{m_ne^{\frac{|\lambda|}{m_n}}}{n^{1/2}}).
\end{aligned}
\end{equation*}
Therefore,
\begin{eqnarray*}
&&\limsup_{n\rightarrow\infty}\frac{m_n^2}{n}\log\mathbf E\Big[\exp\Big\{
\frac{\lambda n}{m_n}\Big(\sum_{i=1}^n\frac 1 n\cos(\sqrt{2nu}\Delta_i^nX)-\int_0^1e^{-u\sigma_s^2}\mathrm{d}s\Big)\Big\}\Big]\\
&\leq&\limsup_{n\rightarrow\infty}\Big\{\lambda m_n\Big(\frac 1 n\sum_{i=1}^ne^{-u\sigma^2_{t_{i-1}}}-\int_0^1e^{-u\sigma_s^2}\mathrm{d}s\Big)
+\frac{\lambda^2}{4n}\sum_{i=1}^n\Big(1+e^{-4u\sigma^2_{t_{i-1}}}\Big)\Big.\\
&&\qquad\quad-\frac{\lambda^2}{2n}\Big(1+\frac{3\lambda^2}{m_n^2}\Big)^{-2}\sum_{i=1}^ne^{-2u\sigma^2_{t_{i-1}}}\Big\}\\
&=&\frac{\lambda^2}{4}\Big(1+\int_0^1e^{-4u\sigma_s^2}\mathrm{d}s-2\int_0^1e^{-2u\sigma_s^2}\mathrm{d}s\Big),
\end{eqnarray*}
where,
\begin{equation}\nonumber
\lim_{m_n\rightarrow\infty}m_n\Big(\frac 1 n\sum_{i=1}^ne^{-u\sigma^2_{t_{i-1}}}-\int_0^1e^{-u\sigma_s^2}\mathrm{d}s\Big)=0,
\end{equation} is derived by the $1/2$-H\"{o}lder continuity of $\sigma_s$.
Similarly, we have
\begin{equation}\nonumber
\begin{aligned}
&\liminf_{n\rightarrow\infty}\frac{m_n^2}{n}\log\mathbf E\Big[\exp\Big\{
\frac{\lambda n}{m_n}\Big(\sum_{i=1}^n\frac 1 n\cos(\sqrt{2nu}\Delta_i^nX)-\int_0^1e^{-u\sigma_s^2}\mathrm{d}s\Big)\Big\}\Big]\\
\geq&
\frac{\lambda^2}{4}\Big(1+\int_0^1e^{-4u\sigma_s^2}\mathrm{d}s-2\int_0^1e^{-2u\sigma_s^2}\mathrm{d}s\Big).
\end{aligned}
\end{equation}
Therefore,
\begin{equation}\label{Lambda1}
\lim_{n\rightarrow\infty}\frac{m_n^2}{n}\log\mathbf E\Big[\exp\Big\{
\frac{\lambda n}{m_n}\Big(\sum_{i=1}^n\frac 1 n\cos(\sqrt{2nu}\Delta_i^nX)-\int_0^1e^{-u\sigma_s^2}\mathrm{d}s\Big)\Big\}\Big]=\Lambda(\lambda,u),
\end{equation}
where $$\Lambda(\lambda,u)=
\frac{\lambda^2}{4}\Big(1+\int_0^1e^{-4u\sigma_s^2}\mathrm{d}s-2\int_0^1e^{-2u\sigma_s^2}\mathrm{d}s\Big).$$
Thus, for $x\in\mathbb R$, $$\Lambda^*(x,u)=\sup_{\lambda\in\mathbb R}\Big\{\lambda x-\Lambda(\lambda,u)\Big\}=\frac{x^2}{1+\int_0^1e^{-4u\sigma_s^2}\mathrm{d}s-2\int_0^1e^{-2u\sigma_s^2}\mathrm{d}s}.$$
Now the result comes from the G\"{a}rtner-Ellis's Theorem and the proof is complete.\hfill $\Box$

\textbf{Proof of Theorem \ref{LDP of process}:}
For $k\in\mathbb N$, $\lambda=(\lambda_1,\cdots,\lambda_k)\in\mathbb R^k$ and
$$\mathcal P=\{0\leq u_1<\cdots<u_k\leq1\},$$
define
\begin{equation*}%\nonumber
%\begin{aligned}
\Lambda_n^{(k)}(\lambda)=\log\mathbf E\Big[\exp\Big\{\sum_{j=1}^k\lambda_jV_n(u_j)\Big\}\Big]
=\log\mathbf E\Big[\exp\Big\{\sum_{i=1}^n\sum_{j=1}^k\frac{\lambda_j}{n}\cos\Big(\sqrt{2u_j n}\Delta_i^nX\Big)\Big\}\Big].
%\end{aligned}
\end{equation*}
Therefore, similar to the proof of Theorem \ref{LDP}, we have
\begin{eqnarray*}
&&\lim_{n\rightarrow\infty}\frac{1}{n}\Lambda_n^{(k)}(\lambda n)\\
&=&\lim_{n\rightarrow\infty}\frac{1}{n}\sum_{i=1}^n\log\mathbf E\Big[\exp\Big\{\sum_{j=1}^k\lambda_j\cos\Big(\sqrt{2u_j n}\Delta_i^nX\Big)\Big\}\Big]\\
&=&\lim_{n\rightarrow\infty}\frac{1}{n}\sum_{i=1}^n\Big[\log\int_{\mathbb R}\frac{1}{\sqrt{2\pi\sigma^2_{t_{i-1}}}}\exp\Big\{-\frac{y^2}{2\sigma^2_{t_{i-1}}}+\sum_{j=1}^k\lambda_j\cos\Big(\sqrt{2u_j}y\Big)\Big\}\mathrm{d}y\Big]\\
&=&\int_0^1\log\int_{\mathbb R}\frac{1}{\sqrt{2\pi\sigma^2_{s}}}\exp\Big\{-\frac{y^2}{2\sigma^2_{s}}+\sum_{j=1}^k\lambda_j\cos\Big(\sqrt{2u_j}y\Big)\Big\}\mathrm{d}y\mathrm{d}s.
\end{eqnarray*}
Let
$$\Lambda^{\mathcal P}(\lambda):=\int_0^1\log\int_{\mathbb R}\frac{1}{\sqrt{2\pi\sigma^2_{s}}}\exp\Big\{-\frac{y^2}{2\sigma^2_{s}}+\sum_{j=1}^k\lambda_j\cos\Big(\sqrt{2u_j}y\Big)\Big\}\mathrm{d}y\mathrm{d}s.$$
Thus, for $k\in\mathbb N$, it follows from G\"{a}rtner-Ellis's Theorem that $\mathbf P((V_n(u_1),\cdots,V_n(u_k))\in\cdot)$ satisfies the LDP with the speed $n$ and the rate function $$I^{\mathcal P}(x)=\sup_{\lambda\in\mathbb R^k}\Big\{\langle\lambda,x\rangle-\Lambda^{\mathcal P}(\lambda)\Big\}.$$
For $r>1$, define $K_r:=[-r,r]$. It is easy to see that $K_r$ is a compact set of $(\mathcal C([0,1];\mathbb R),\|\cdot\|)$.
For $\lambda>0$,
\begin{equation}\nonumber
\begin{aligned}
\frac{1}{n}\log\mathbf P\Big(\sup_{0\leq u\leq1}|V_n(u)|>r\Big)\leq\frac{1}{n}\log\frac{\mathbf E\Big[\exp\Big\{\sup_{0\leq u\leq1}n\lambda|V_n(u)|\Big\}\Big]}{e^{n\lambda r}}\leq-\lambda(r-1).
\end{aligned}
\end{equation}
Therefore, $\mathbf P(V_n(\cdot)\in\cdot)$ is exponentially tight on $(\mathcal C([0,1];\mathbb R),\|\cdot\|)$.
Thus, by Proposition I.5.5 of \citet{W1997}, $\mathbf P(V_n(\cdot)\in\cdot)$ satisfies the LDP on $(\mathcal C([0,1];\mathbb R),\|\cdot\|)$ with the speed $1/n$ and the rate function given by
$$I^\infty(\phi)=\sup_{\mathcal P}I^{\mathcal P}(\phi(\mathcal P)),\quad\phi\in\mathcal C([0,1];\mathbb R),$$
where the supremum is taken over all finite partitions $\mathcal P=\{0\leq u_1<\cdots<u_k\leq1\}$ of $[0,1]$. For fixed $k\in\mathbb N$, let $\mathcal P^k$ denote the finite partition $0=u_0<u_1<\cdots<u_k=1$ such that $u_j=\frac{j}{k}$ for $j=0,\cdots,k$.
Let $Y\sim N(0,\sigma_s^2)$, it follows from H\"{o}lder inequality that
\begin{eqnarray*}
I^{\mathcal P^k}(\phi(\mathcal P^k))&=&\sup_{\lambda\in\mathbb R^k}\Big\{\sum_{j=1}^k\lambda_j\phi(u_j)-\int_0^1\log\int_{\mathbb R}\frac{1}{\sqrt{2\pi \sigma_s^2}}\exp\Big\{-\frac{y^2}{2\sigma_s^2}+\sum_{j=1}^k\lambda_j\cos(\sqrt{2u_j}y)\Big\}\mathrm{d}y\mathrm{d}s\Big\}\\
&=&\sup_{\lambda\in\mathbb R^k}\Big\{\sum_{j=1}^k\lambda_j\phi(u_j)-\int_0^1\log\mathbf E\Big[\prod_{j=1}^k\exp\Big\{\lambda_j\cos(\sqrt{2u_j}Y)\Big\}\Big]\mathrm{d}s\Big\}\\
&\geq&\sup_{\lambda\in\mathbb R^k}\Big\{\sum_{j=1}^k\lambda_j\phi(u_j)-\int_0^1\log\prod_{j=1}^kE\Big[\exp\Big\{k\lambda_j\cos(\sqrt{2u_j}Y)\Big\}\Big]^{1/k}\mathrm{d}s\Big\}\\
&=&\sup_{\lambda\in\mathbb R^k}\Big\{\sum_{j=1}^k\Big(\lambda_j\phi(u_j)-\frac{1}{k}\Lambda(k\lambda_j,u_j)\Big)\Big\}\\
&=&\sum_{j=1}^k\frac{1}{k}\sup_{\lambda_j\in\mathbb R}\Big(k\lambda_j\phi(u_j)-\Lambda(k\lambda_j,u_j)\Big)\\
&=&\sum_{j=1}^k\frac{1}{k}I(\phi(u_j),u_j).
\end{eqnarray*}
where $\Lambda(\cdot,\cdot)$ and $I(\cdot,\cdot)$ are defined in \eqref{limiting logarithmic generating function} and \eqref{rate function}.
Therefore, $$\hspace{5cm} I^\infty(\phi)\geq\int_0^1I(\phi(u),u)du.\hspace{5cm} \Box$$

\

\textbf{Proof of Theorem \ref{MDP of process}:}
For $k\in\mathbb N$, $\lambda=(\lambda_1,\cdots,\lambda_k)\in\mathbb R^k$ and $\mathcal P=\{0\leq u_1<\cdots<u_k\leq1\}$, similar as Theorem \ref{MDP} and Theorem \ref{LDP of process}, we have
\begin{equation}\nonumber
\begin{aligned}
&\lim_{n\rightarrow\infty}\frac{m_n^2}{n}
\log\mathbf E\Big[\exp\Big\{\sum_{j=1}^k\frac{\lambda_jn}{m_n}\Big(\sum_{i=1}^n\frac{1}{n}\cos\Big(\sqrt{2u_j n}\Delta_i^nX\Big)
-\int_0^1e^{-u_j\sigma_s^2}\mathrm{d}s\Big)\Big\}\Big]\\
=&\frac{1}{4}\sum_{j=1}^k\lambda_j^2\int_0^1\Big(1-e^{-2u_j\sigma_s^2}\Big)^2\mathrm{d}s
+\frac{1}{2}\sum_{1\leq j<l\leq k}\lambda_j\lambda_l\int_0^1\Big(e^{-\frac{(\sqrt{u_j}+\sqrt{u_l})^2\sigma_s^2}{2}}
-e^{-\frac{(\sqrt{u_j}-\sqrt{u_l})^2\sigma_s^2}{2}}\Big)^2\mathrm{d}s\\
=:&\Lambda^\mathcal P(\lambda).
\end{aligned}
\end{equation}
Thus, for $k\in\mathbb N$, it follows from G\"{a}rtner-Ellis's Theorem that $$\mathbf P\Big(\Big(m_n\big(V_n(u_1)-F(u_1)\big),\cdots,m_n\big(V_n(u_k)-F(u_k)\big)\Big)\in\cdot\Big)$$ satisfies the LDP with the speed $m_n^2/n$ and the rate function $$I^{\mathcal P}(x)=\sup_{\lambda\in\mathbb R^k}\Big\{\langle\lambda,x\rangle-\Lambda^{\mathcal P}(\lambda)\Big\}.$$
It follows from Corollary I.5.3 of \citet{W1997} that $\mathbf P(m_n(V_n(\cdot)-F(\cdot))\in\cdot)$ satisfies the weak LDP on $\mathcal C([0,1];\mathbb R)$ with respect to the pointwise convergence topology with the speed $m_n^2/n$ and the rate function given by
$$I^\infty(\phi)=\sup_{\mathcal P}I^{\mathcal P}(\phi(\mathcal P)),\quad\phi\in\mathcal C([0,1];\mathbb R),$$
where the supremum is taken over all finite partitions $\mathcal P=\{0\leq u_1<\cdots<u_k\leq1\}$ of $[0,1]$. \qed

\

\textbf{Proof of Theorem \ref{LDP2}:} The proof is similar to that of Theorem \ref{LDP}, but now we need to compute
\begin{eqnarray*}
& &\lim_{n\rightarrow\infty}\frac{1}{n}\log\mathbf E\Big[\exp\Big\{n\lambda\sum_{i=1}^{N_n}\Delta_i^nt\cos(\sqrt{2u/\Delta_i^nt}\Delta_i^nX)\Big\}\Big]\\
&=&\lim_{n\rightarrow\infty}\sum_{i=1}^{N_n}\log\mathbf E\Big[\exp\Big\{\lambda n\Delta_i^nt\cos(\sqrt{2u/\Delta_i^nt}\Delta_i^nX)\Big\}\Big]\frac{\Delta_i^nt}{n\Delta_i^nt}.
\end{eqnarray*}
Noting that $(n\Delta_i^nt)^{-1}=T_n'(t_{n,i-1}+)$, hence the result in 1) follows from the proof of Theorem \ref{LDP} and Assumption \ref{asu:time2}.

To prove the result in 2), by repeating the proof of Theorem \ref{MDP}, but now instead of \eqref{Lambda1}, we have
\begin{equation}
\lim_{n\rightarrow\infty}\frac{m_n^2}{n}\log\mathbf E\Big[\exp\Big\{
\frac{\lambda n}{m_n}\Big(\sum_{i=1}^{N_n}\Delta_i^nt\cos(\sqrt{2u/\Delta_i^nt}\Delta_i^nX)-\int_0^1e^{-u\sigma_s^2}\mathrm{d}s\Big)\Big\}\Big]=\Lambda(\lambda,u),
\end{equation}
where $$\Lambda(\lambda,u)=
\frac{\lambda^2}{4}\int_0^1 \frac{1+e^{-4u\sigma_s^2}-2e^{-2u\sigma_s^2}}{T'(s)}\mathrm{d}s.$$
This completes the proof of 2).
\hfill $\Box$

\

\textbf{Proof of Theorem \ref{LDP3}}: The proof is simply a repeat of proof of Theorems \ref{LDP of process} and \ref{MDP of process}, just with some changes of the notation. \hfill $\Box$

\section*{Acknowledgement}
The authors would like to thank the editor, an associate editor, and an anonymous referee for their very extensive and
constructive suggestions that helped to improve this paper considerably. Liu gratefully acknowledges financial support from The Science and Technology Development Fund, Macau SAR (No. 202/2017/A3) and NSFC (No. 11971507).

\section*{References}
\bibliographystyle{model2-names}
\bibliography{liu}
\end{document}